\documentclass{amsart}
\usepackage{amssymb}
\usepackage{mathdots}
\usepackage[...]{youngtab}
\usepackage{tikz}
\usetikzlibrary{shapes,snakes}
\usepackage{framed}
\usepackage{tikz-cd}
\usepackage{amsmath}
\usepackage{txfonts}
\usepackage{adjustbox}
\usepackage{booktabs}

\input xy
\xyoption{all}

\usepackage{amsfonts}
\usepackage{mathrsfs}

\usepackage{latexsym}
\usepackage{graphicx}
\usepackage{amscd,amssymb,amsmath,amsbsy,amsthm}
\usepackage[all]{xy}

\usepackage{xcolor}
\definecolor{darkred}{RGB}{139,0,0}
\usepackage{hyperref}
\hypersetup{colorlinks,linkcolor={darkred},citecolor={darkred},urlcolor={darkred}} 
\usepackage{verbatim}
\usepackage{tikz-cd}
\usepackage{amsmath}
\usepackage{amsfonts}
\usepackage{amssymb}
\usepackage{pdfpages}
\usepackage{tikz}
\usepackage{asymptote}
\usetikzlibrary{decorations.markings,arrows}
\usepackage[all]{xy}
\usepackage{yfonts}
\usepackage[...]{youngtab}

\usepackage{mathrsfs}
\usepackage{latexsym}
\usepackage{graphicx}
\usepackage{amscd,amssymb,amsmath,amsbsy,amsthm}
\usepackage[all]{xy}
\newtheorem{theorem}{Theorem}[section]

\theoremstyle{definition}
\newtheorem{defn}[theorem]{Definition}
\newtheorem{ex}[theorem]{Example}
\newtheorem{remark}[theorem]{Remark}
\newtheorem{recall}[theorem]{Recall}

\newtheorem{conj}[theorem]{Conjecture}

\newsavebox{\boxedtikzcdbox}
\DeclareMathOperator{\Gr}{Gr}
\DeclareMathOperator{\id}{\mathbf{id}}

\DeclareMathOperator{\FF}{\mathbb{F}}
\DeclareMathOperator{\NN}{\mathbb{N}}
\DeclareMathOperator{\ZZ}{\mathbb{Z}}
\DeclareMathOperator{\QQ}{\mathbb{Q}}

\DeclareMathOperator{\CC}{\mathbb{C}}
\DeclareMathOperator{\PP}{\mathbb{P}}

\DeclareMathOperator{\DD}{\mathbb{D}}

\DeclareMathOperator{\rad}{\mathbf{rad}}
\DeclareMathOperator{\cG}{\mathcal{G}}

\DeclareMathOperator{\cN}{\mathcal{N}}

\DeclareMathOperator{\pO}{\hat{\mathcal{O}}}
\DeclareMathOperator{\hO}{\widehat{\mathcal{O}}}
\DeclareMathOperator{\cK}{\mathcal{K}}

\DeclareMathOperator{\Fl}{\mathbf{Fl}}

\DeclareMathOperator{\Mod}{\mathbf{Mod}}

\DeclareMathOperator{\Span}{\mathbf{span}}

\DeclareMathOperator{\cO}{\mathcal{O}}

\DeclareMathOperator{\ttop}{\mathbf{top}}

\DeclareMathOperator{\fm}{\mathfrak{m}}
\DeclareMathOperator{\fF}{\mathfrak{F}}

\DeclareMathOperator{\fb}{\mathfrak{b}}
\DeclareMathOperator{\fn}{\mathfrak{n}}

\DeclareMathOperator{\fh}{\mathfrak{h}}

\DeclareMathOperator{\Sym}{\mathbf{Sym}}

\DeclareMathOperator{\bd}{\mathbf{d}}

\DeclareMathOperator{\fG}{\mathfrak{G}}
\DeclareMathOperator{\SL}{\mathbf{SL}}
\DeclareMathOperator{\GL}{\mathbf{GL}}
\DeclareMathOperator{\sln}{\mathfrak{sl}}
\DeclareMathOperator{\gl}{\mathfrak{gl}}
\DeclareMathOperator{\fg}{\mathfrak{g}}
\DeclareMathOperator{\Spec}{\mathbf{Spec}}
\DeclareMathOperator{\rank}{\mathbf{rank}}

\DeclareMathOperator{\ord}{\mathbf{ord}}
\DeclareMathOperator{\Mat}{\mathbf{Mat}}
\DeclareMathOperator{\Hom}{\mathbf{Hom}}

\DeclareMathOperator{\sym}{\mathcal{\mathbf{Sym}}}
\DeclareMathOperator{\cF}{\mathcal{\mathbf{F}}}

\DeclareMathOperator{\cA}{\mathbf{\mathcal{A}}}

\newcommand{\la}{\langle}
\newcommand{\ra}{\rangle}

\newcount\cols
{\catcode`,=\active\catcode`|=\active
	\gdef\Young(#1){\hbox{$\vcenter
			{\mathcode`,="8000\mathcode`|="8000
				\def,{\global\advance\cols by 1 &}%
				\def|{\cr
					\multispan{\the\cols}\hrulefill\cr
					&\global\cols=2 }%
				\offinterlineskip\everycr{}\tabskip=2pt
				\dimen0=\ht\strutbox \advance\dimen0 by \dp\strutbox
				\halign
				{\vrule height \ht\strutbox depth \dp\strutbox##
					&&\hbox to \dimen0{\hss$##$\hss}\vrule\cr
					\noalign{\hrule}&\global\cols=2 #1\crcr
					\multispan{\the\cols}\hrulefill\cr%
				}
			}$}}
}

\title[Surface Algebras and Surface Orders II: Affine Bundles on Curves]
{Surface Algebras and Surface Orders II: Affine Bundles on Curves}

\author[A.~Schreiber]{Amelie~Schreiber}
\email{amelie.schreiber.math@gmail.com}

\subjclass[2010]{
	Primary
	11G32  	%Dessins d'enfants, Belyĭ theory
	16E05  	%Syzygies, resolutions, complexes
	Secondary
	05E10  %Combinatorial aspects of representation theory
	05E15  	%Combinatorial aspects of groups and algebras
}

\date{\today}
\keywords{Riemann surfaces, surface algebra, surface order, loop group, loop algebra, Langlands Correspondence, affine Grassmannian, affine flag variety}

\begin{document}

	\begin{abstract}
In this paper I provide a way to understand the Geometric Langlands Program. 
	\end{abstract}

\maketitle

\textcolor{darkred}{\hrule}

\tableofcontents

\textcolor{darkred}{\hrule}
\medskip
\section{Introduction}

\subsection{The Big Picture}
In this paper, I will present a construction of a gluing of loop algebras and loop groups given by a combinatorial object associated to a covering of the projective line $\PP_{\CC}^1$ by some compact Riemann surface $X$, which I called a \emph{constellation} in my last paper \cite{S1}. In particular, there will be a loop group associated to a punctured disk $\DD^{\times}_{x_j} = \Spec \CC((x_j))$, for each ramification point of $X$. At times, we will call the $x_j$ ramification points, even though we will need unramified topological coverings at times, in which case they are to be thought of as punctures of the surface $X$. I will sometimes refer to this as the "unramified case."\footnote{I will use this terminology, even though this could just mean we are at a point $y$ in the base space with a unique point in the fibre of the covering, but we are not really interested in such points at the moment.} Once this association of loop groups and loop algebras is constructed, I will explain how one can define a fibre bundle over $X$ which has an embedding into a trivial bundle of a product of loop groups. Once this bundle structure is explained, we can then look at representations of the loop algebras of a covering of Riemann surfaces and their associated pullback which defines a surface order. I will then show that the representations of the surface orders can be identified with affine Schubert varieties inside the pullback of affine flag varieties. 

It will be shown that these spaces luckily have a very nice structure which can be well explained in terms of the affine Schubert varieties. Moreover, we have that the pullback of loop algebras describes a pullback of the spaces $C^{\infty}(S^1, G((x_i)))$, where $S^1$ is the circle and $G((x_i))$ are the loop groups of interest. This, somewhat unsurprisingly, gives a map of the monodromy group of the covering into the pullback. We may think of this as a gluing of copies of $S^1$, each copy being a deformation retract of the punctured discs associated to one of the ramification points/punctures $x_j$, and identified with an element of the fundamental group(oid) of the surface. It can also be thought of as a gluing of cyclic groups, each corresponding to a cyclic subgroup of the monodromy group of the covering. Given any covering $p: X \to Y$, of compact Riemann surfaces, this allows us to interpret the Galois group of the field extension given by $\CC(Y) \to \CC(X)$, in terms of the representation theory and geometry of loop groups and loop algebras. In particular, take the affine Kac-Moody Lie algebras $\widehat{\fg}_j$ given by central extensions
	\[ \CC \mathbf{1} \to \widehat{\fg}_j \to \fg((x_j)) \] 
of each loop algebra $\fg((x_j))$. Next, the pullback of the various loop algebras $\fg((x_j))$ and their Borel subalgebras $\fb((x_j))_{+}$, each associated to a punctured disk $\DD_j^{\times} \subset X$ over the (would be) ramification point $x_i$ (i.e. punctures for an unramified covering), gives us a way of understanding the Galois groups $\cG(\overline{\FF}_j/\FF_j) \cong \hat{\ZZ}$, where $\FF_j = \CC((x_j))$. To be precise, we can associate an automorphism of $\FF((x_j^{1/n}))$, defined by 
\[ x_j \mapsto e^{2 \pi i k/n}\cdot x_j \]
to an element $\sigma_j \in \cG\left(\FF((x_j^{1/n_j}))/\CC((x_j))\right) \cong \ZZ/n_j\ZZ \hookrightarrow S_{N}$, where $n_j$ is the degree of the ramification at $x_j$.\footnote{The positive integer $n_j$ is also the size of the hereditary orders of the normalization of the surface order for $X$, and the number of arrows in the cyclic quiver of the normalization of the surface algebra of $X$} This corresponds to a cyclic subgroup of the monodromy group of $X$, embedded as a permutation group in a symmetric group, and given as the pullback of the cyclic groups $\la \sigma_j \ra \cong \ZZ/n_j\ZZ$, over each ramification point/puncture $x_j$. 

This gives us a way of understanding the monodromy group of the Riemann surface $X$ as a pullback of cyclic subgroups coming from automorphisms of $\FF((x_j^{1/n}))$, where the algebraic completion of $\FF_j$ is given by the inductive limit $\varprojlim \FF_j((x_j^{1/n}))$, of the directed system given by inclusions
\[ \FF_j((x_j^{1/n})) \hookrightarrow \FF_j((x_j^{1/m})) \]
if $n$ divides $m$. More generally, for any covering $p:X \to Y$ of Riemann surfaces, we have a map of fields of regular functions
\[ \CC(Y) \to \CC(X) \]
and $\cG\left(\CC(X)/\CC(Y)\right)$ may be identified with the group of Deck transformations. For an \emph{unramified} cover, we can identify this with a quotient of $\pi_1(Y)$. If we are in the situation where 
\[ Y = \PP^1_{\CC} -\{0,1\} \]
is the Riemann sphere with two punctures, then $\pi_1(Y) = \fF_2 = \ZZ * \ZZ = \la x, y \ra$, is a free group on two generators. The fascinating fact that will come from all of this is the following:\\
\\
The universal covering of every surface algebra is the bi-colored, directed, Cayley graph of $\fF_2$. This means we may identify completions of surface algebras, i.e. the pullbacks of Borel subalgebras $\fb((x_i))$ of loop algebras $\fg((x_i))$, over punctured disks $\DD_j^{\times}$ at ramification points (or punctures) $x_j$, with finite index quotients of the free group on two generators given by coverings 
\[ p: X \to \PP^1_{\CC}-\{0,1\}. \]

\subsection{Layout of the Paper}
In the first part of the paper, Section \ref{Gel'fand Ponomarev} I will give a sketchy but (hopefully) very intuitive exposition of some of the results of \cite{GP} by Gel'fand and Ponomarev, and this gives us a first pass at the big ideas, with some concrete examples we can get our hands on. The \cite{GP} paper, in many ways, led to so many ideas in the representation theory of quivers, the representation theory of Lie algebras, and later to Kac-Moody Lie algebras. Its influence can be seen in physics and the study of conformal field theories, and there are likely many other places where this paper sparked ideas that are adrift in the aether for the time being. The section on the \emph{Gel'fand-Ponomarev Surface Algebra} contains the most basic examples of the ideas that we will need for much of the rest of the paper. The combinatorics, intuition, and the way of thinking about this example can be seen hiding in the background of countless papers following Gel'fand and Ponomarev's example. Many of these ideas are obscured, never explicitely stated, or have been delivered in such high brow and exclusive language as to be inaccessible to anyone not quite adept at or steeped in the culture of the representation theory of quivers. For those looking to gain access to these ideas, which are fairly ubiquitous in the quiver community, this section alone may be a good read to develop intuition and see examples not present in the literature, unless of course you know David Benson.\footnote{David Benson's infamous $60+$ GB file, colloquially known as "\emph{the Benson Archive}", is a mass collection of math papers and books, some of which are quite rare and old, and many seem to exists nowhere else on the planet, for example \cite{G}.}

Once I have "\emph{explained}" the Gel'fand-Ponomarev Surface algebra (which, quite importantly, is the surface algebra of the trivial dessin d'enfant $\PP_{\CC}^1-\{0,1\}$) I will set the conventions for the Lie groups, Lie algebras, and the associated loop groups and loop algebras used. Then we can take a look at some of the main ideas we are going to try to understand here. Then I will provide a brief review of the \emph{Surface Algebras} and \emph{Surface Orders} that were defined in my previous paper \cite{S1} \emph{"Surface Algebras and Surface Orders I: Homological Characterizations of Dessins D'Enfants and Riemann Surfaces"}. It will be necessary to understand that material to understand this paper. Next, I will give a brief review of Affine Grassmannian, Affine Flag Varieties, and Affine Schubert Varieties. We follow mostly Magyar's notes \cite{M}, and also Bj\"{o}rner and Brenti's book \emph{"The Combinatorics of Coxeter Groups"} \cite{BB}. At this point we can look at the basic objects of interest, which I will call \emph{Trivial Schubert Bundles}. Next, we will look at a fibre product of Affine Schubert Varieties which will give us a nontrivial bundle over the complex curve $X$ (or Riemann surface). Once we have define this bundle, I will show that its geometry can actually be well understood from a different perspective, as well as from the perspective of Affine Schubert Varieties. This setup will be necessary for the next paper in which we describe some moduli spaces of representations of the surface algebras and surface orders, which are related to the Affine Schubert Bundles and are incredibly rich in structure. The Geometric Invariant Theory (GIT) and representation theory needed for this is simply too much to include here, and requires its own lengthy setup, so we save that for another day. 

\emph{Currently, there are no "Theorems". This is intentional and I intend keep things this way as it removes emphasis from trying to prove important results, and puts it squarely on writing an interesting, useful, and accessible paper. It also feels much less natural to break things up into blocks. I prefer there to be a kind of flow to my work. It is likely many will frown upon this. I don't care. I want my work to be genuine and in my mind, each paper is actually one long "Theorem", chiseled off of the bigger picture I am hoping to paint. I'm following the math, not conventions, and I intend to be true to the math.}

\subsection{A Quick Review of Constellations, Monodromy, and Ramified Coverings}
This short aside is merely to set the terminology and conventions used throught the paper. So, 

Let $S_{2N}$ be the symmetric group on $[2N] = \{1,2,3...,2N\}$. Permutations will act on the left, and if $\sigma \in S_{2N}$, we will use the notation $\sigma \cdot i = \sigma(i)$. For example, for $\sigma = (1,3,2) \in S_3$ we have
\[ \sigma(1) = 3, \quad \sigma(2) = 1, \quad \sigma(3) = 2.\]
Let us define a $k$-\textbf{constellation} to be a sequence $C=[g_1, g_2, ..., g_k]$, $g_i \in S_n$, such that:
\begin{enumerate}
\item The group $G = \la g_1, g_2, ..., g_k \ra$ generated by the $g_i$ acts \textit{transitively} on $[n]$. 
\item The product $\prod_i^k g_i = \id$ is the identity. 
\end{enumerate}
The constellation $C$ has "\textbf{degree} $n$" in this case, and "\textbf{length} $k$". Our main interest will be in $3$-constellations $C = [\sigma, \alpha, \phi]$, with $\alpha$ a fixed-point free involution. The group $G = \la \sigma, \alpha \ra = \la \sigma, \alpha, \phi \ra$ since $\sigma \alpha = \phi^{-1}$, and $\alpha$ being an involution means $\alpha^{-1} = \alpha$. This group will be called the \textbf{cartographic group} or the \textbf{monodromy group} generated by $C$. 

Now, let $p:X^{\circ} \to Y$ be a covering of a compact Riemann surface $Y$ by some topological surface $X^{\circ}$. It is a well known fact that is there is a finite sheeted covering
\[ f^{\circ}: X^{\circ} \to Y-\{y_1, ..., y_t\} \]
then up to isomorphism there is a unique compact Riemann surface $X \supset X^{\circ}$, so that $X^{\circ}$ is dense, and in fact is just $X$ with a finite number of punctures $\{x_1, ..., x_r\}$, such that $f^{\circ}$ extends to a holomorphic map of Riemann surfaces
\[ f:X \to Y \]
Generally speaking the set $B_Y = \{y_1, ..., y_t\}$ is the branch locus of $Y$, and every $x_i$ lies in the fibre over some $y_j$. Let us call the extended holomorphic map $f:X \to Y$ a \textbf{ramified cover} of $Y$. Let us call the map $f^{\circ}: X^{\circ} \to Y-\{y_1, ..., y_t\}$ an \textbf{unramified cover}. So, locally, neighborhoods of the $x_ \in X^{\circ}i$ are open disks with punctures at $x_i$. We will often blur this distinction throughout, and call the point which is added to $X^{\circ}$ to obtain the compact surface $X$ also $x_i$. It should be clear from the context which case is relevant and would become quite cumbersome to constantly distinguish between the two throughout the entire paper. For example, when we speak of punctured disks $\DD_i^{\times}$ around a \textbf{ramification point} $x_i$, we mean around a puncture in the unramified cover 
\[ f^{\circ}: X^{\circ} \to Y-\{y_1, ..., y_t\}. \]

\emph{For more information on the combinatorics of such things, how they define cellularly embedded graphs on Riemann surfaces, and various homological characterizations of them please see \cite{S1}.}

\section{Combinatorics of the Gel'Fand Ponomarev Algebra}\label{Gel'fand Ponomarev}

\subsection{The Gel'fand-Ponomarev Algebra}
In the now classic paper \cite{GP}, Gel'fand and Ponomarev studied the indecomposable representations of the Lorentz group, which is equivalent to classifying the "\emph{Harish-Chandra modules}" of the Lie algebra $\sln_2(\CC)$. 
\footnote{For an explanation of the finite dimensional representations of this Lie algebra and examples, I highly recommend \cite{FH} Chapter 10 and 11.}
They used methods now standard in the representation theory of so-called \emph{string algebras} and more generally \emph{special biserial algebras}.\footnote{For an excellent introduction to these algebras which will provide more insite into how this relates to modular representation theory of finite groups I recommend \cite{E} and \cite{B1, B2}.}
In their study of Harish-Chandra modules, they studied the following question:

\emph{Let $k$, be a field and let $V_1$ and $V_2$ be finite dimensional $k$-vector spaces. Suppose we have \textbf{nilpotent operators}
\[ X : V_1 \to V_2, \quad Y: V_2 \to V_1, \quad H: V_2 \to V_2 \]
such that $YH = 0 = HX$. Fix a basis of $V_1$ and $V_2$ and classify all canonical forms of $X$ and $Y$ and $H$. 
}
In \cite{FH} one sees the following diagram for $\sln_2(\CC)$
\[ \xymatrix{
\cdots \ar[r]_{X} & V_{\alpha-4} \ar@(ul,ur)^{H} \ar[r]_{X} \ar@/_/[l]_{Y} & V_{\alpha-2} \ar@(ul,ur)^{H} \ar[r]_{X} \ar@/_/[l]_{Y} & V_{\alpha} \ar[r]_{X} \ar@/_/[l]_{Y} \ar@(ul,ur)^{H}  & V_{\alpha+2} \ar@(ul,ur)^{H} \ar[r]_{X} \ar@/_/[l]_{Y} & \ar@/_/[l]_{Y} \cdots
}\]
with the typical basis
\[ X = \begin{pmatrix}
 0 & 1 \\
 0 & 0
 \end{pmatrix}, \quad Y=\begin{pmatrix}
 0 & 0 \\
 1 & 0
 \end{pmatrix}, \quad H = \begin{pmatrix}
 1 & 0 \\
 0 & -1
 \end{pmatrix}\]
for $\sln_2(\CC)$. Now, the idea, is we wish to send $\sln_2(\CC) \hookrightarrow \sln_n(\CC)$ in a meaningful way and study the action. If we take instead the algebra generated by
 \[ X = \begin{pmatrix}
 0 & Y \\
 X & 0
 \end{pmatrix}, \quad Y = \begin{pmatrix}
 0 & H \\
 0 & 0
 \end{pmatrix} \]
 acting on 
\[V = V_1 \oplus V_2\]
then this is equivalent to classifying all representations of the following path algebra:
$I=\la xy, yx \ra$, $\Lambda = kQ/I$, 

\[ \xymatrix{
	\bullet \ar@(ul,dl)_{x} \ar@(dr,ur)_{y}
} \]
What is remarkable is that we can recover the action of $\sln_2(\CC)$ via the embedding into $\sln_n(\CC)$, using MacLane's "linear relations", as explained by Gel'fand and Ponomarev.\footnote{A more thorough explanation of this is given by Ringel in \cite{R1} on his study of the representations of the infinite dihedral group $D_{\infty}$ and its quotients. This is more or less the same problem, so we will get into that another time, when it proves useful.}

The \textbf{path algebra} of the Gel'fand-Ponomarev algebra is then 
\[ k\la  x, y \ra  /\la xy, yx \ra\cong k[x,y]/(xy).\] 
One may then take the \textbf{dessin order} \footnote{The path algebra and dessin and surface orders are defined in \cite{S1} and in the forthcoming sections of this paper.}
to be $\Lambda = R_{1,2} = k[x,y]/(xy)$ which has normalization $R_1 \times R_2 = k[x] \times k[y]$, with maximal ideals $\fm_1 = (x)$ and $\fm_2 = (y)$ respectively. Fixing an isomorphism $R_1/\fm_1 = k \to k = R_2/\fm_2$, we get a gluing
\[ \xymatrix{
R_{1,2} = k[x,y]/(xy) \ar[r] \ar[d] & k[y] \ar[d] \\
k[x]  \ar[r] & k
}\]
The completion of the path algebra $\Lambda = kQ/I$ is then $\widehat{\Lambda} = k[[x,y]]/(xy)$. One of the main ideas of this paper is not just to generalize what Gel'fand and Ponomarev did, but to uncover deeper number theoretic, geometric, and group theoretic consequences as well. 

\subsection{Bi-colored, Directed Cayley Graph of $\fF_2$}

The directed bi-colored Cayley graph of the free group on two generators 
\[ \fF_2 = \ZZ * \ZZ = \la x, y \ra \]
\emph{\textbf{is actually the universal cover of all surface algebras}}\footnote{We will come back to a proof of this an apply it to surface algebras in general to get our hands on a better understanding of coverings of Riemann surfaces, the corresponding field extensions, how this relates to coverings of algebras and monodromy groups, and how this all gives information in the loop groups and loop algebras in Section \ref{Coverings of Surface Algebras}.} and if, we take the relations in the corresponding infinite path algebra of the Cayley graph to be
\[ k \la x,y, x^{-1}, y^{-1} \ra /\la xy, yx \ra . \]
This is just a quotient of the \textbf{group algebra}, $k\fF_2$ of $\fF_2$. Note, we have an exact sequence of algebras\footnote{This is actually quite important and has implications to the modular representation theory of finite groups and their "blocks", i.e. indecomposable ideals of the group algebra $RG$ as $RG$-modules, for some commutative ring $R$. In particular, if we work in the case of $\QQ_{(p)}, \ZZ_{(p)}, \ZZ/p$, and study quotients of this group algebra by actions of fundamental groups of Riemann surfaces covering the sphere $\PP^1_{\QQ_{(p)}}$, we will obtain "Brauer graph algebras", each of which is a quotient of a surface algebra. Additionally, we can get into the theory of lattices over orders as explained in \cite{CR1, CR2} applied to modular representations of group algebras.}
\[ k \la x,y^{-1} \ra \amalg k \la x^{-1}, y \ra \to k\fF_2 \to k\fF_2/\la xy, yx \ra \]

In what follows, abuse of notation may occur, and we will take $k \fF_2 /I$ to mean the (free) path algebra of the directed Cayley graph of $\fF_2$, modulo some two sided ideal $I$. For a group acting on the graph (or the path algebra), we should also think of $\fF_2/\fG$ as the quotient by the group action on the graph.

Let us now represent the letter "$x$" by a \emph{blue} arrow, and by a \emph{red} arrow we will represent the letter "$y$". Picking an arbitrary basepoint, call it $e_{0}$, in the bi-colored, four-regular, directed tree (see Figure \ref{Cantor Space}), we have arbitrary  words in the alphabet $\{x, y^{-1}\} \amalg \{x^{-1}, y\}$ which can be formed. The choice of the basepoint for $e_{0}$ is of course arbitrary since this graph is a model of the self-similar fractal Cantor space. In particular, we may identify this graph, which with some mild abuse of notation we will call $\fF_2$, with $2^{\omega} \amalg 2^{\omega} = \{0,1\}^{\omega} \amalg \{0,1\}^{\omega}$. Indeed, initially, at our chosen base point $e_{0}$, we have \emph{four} choices of which "direction" to travel along the quiver $\fF_2$. We can travel in one of the four directions "$x, y, x^{-1}, y^{-1}$". Once we have chosen the initial direction, we cannot reverse our direction along the arrow just traversed. In other words, we allow only "\emph{reduced words}". So, once we have chosen to travel either to the right or to the left, we must continue to do so, and our only choices at that point are to go either "up" or "down". 

\medskip

\begin{center}\label{Cantor Space}
\fbox{\includegraphics[
    page=1,
    width=200pt,
    height=200pt,
    keepaspectratio
]{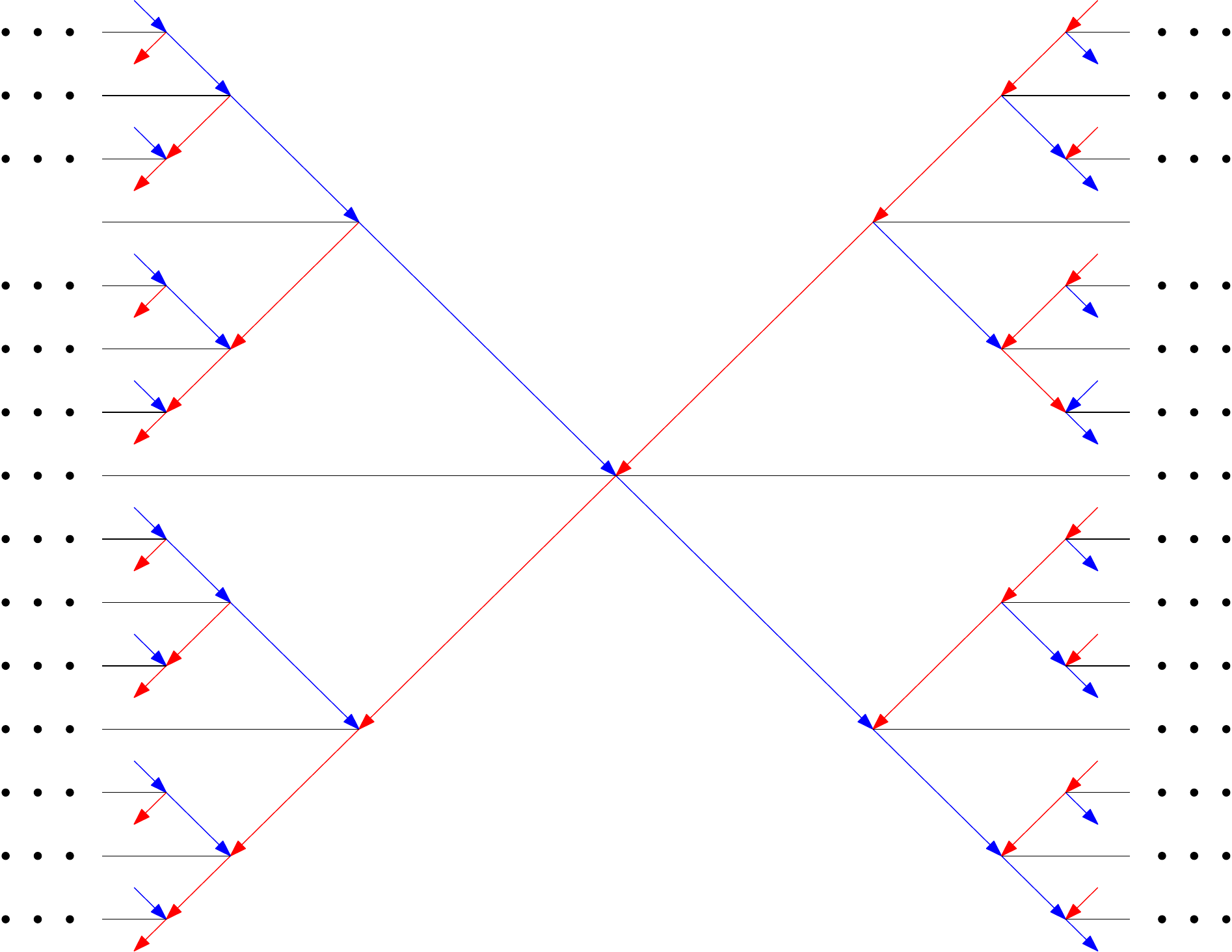}}
\end{center}

\medskip

Now, let us define a linear order on all words. Once a base point $e_0$ is chosen, we will define 
\[  y < x < e_0 < x^{-1} < y^{-1}. \] 
Next, for any word $w$ traveling to the right in the quiver $\fF_2$, i.e. of the form $\ell \cdot y^{-1}$ or $\ell \cdot x$, we will say $x \cdot w < w < y^{-1} \cdot w$. For any word traveling to the left in the quiver $\fF_2$ we will define $ y \cdot w < w < x^{-1} \cdot w$. In less formal terms this means extending to the right is always "\emph{larger than}" extending to the left. Similarly, extending the word upward is always "\emph{larger than}" extending downward. In terms of binary sequences in Cantor space this is just the obvious lexicographic order. 

\begin{ex}
Let $w = x^2y^{-3}x^3y^{-2}$. Then the binary sequence would be
\[  s(w) = 1100011100 \hookrightarrow \{0,1\}^{\omega} \amalg \{0,1\}^{\omega} \]
This is a finite binary string which can be embedded into the disjoint union of two copies of Cantor space via the map
\[ s(w) \mapsto \emptyset \amalg \{1100011100 \cdot 000\cdots\} \subset \{0,1\}^{\omega}  \amalg\{0,1\}^{\omega}.\]
\end{ex}

\subsection{What Does Combinatorial Commutative Algebra Have to Say About This?}

One can think of this as a kind of splitting of the modules over $\CC[[X,Y]]/(XY)$. In combinatorial commutative algebra, it is useful to view modules over $\CC[X,Y]$ as subsets in the lattice $\NN^2$ of $\ZZ^2$. For example, the free $\CC[X,Y]$-module $(X Y^2)$ would be represented graphically as in Figure \ref{Free C[X,Y]-module}. 

\medskip
\adjustbox{scale=0.7,center}{
\begin{tikzcd}\label{Free C[X,Y]-module}
\vdots & \vdots & \vdots & \vdots &  \\
x^3 \arrow[u] \arrow[r] & x^3y \arrow[r] \arrow[u] & \fbox{$x^3y^2$} \arrow[r, dotted] \arrow[u, dotted] & \fbox{$x^3y^3$} \arrow[r, dotted] \arrow[u, dotted] & \cdots \\
x^2 \arrow[u] \arrow[r] & x^2y \arrow[r] \arrow[u] & \fbox{$x^2y^2$} \arrow[r, dotted] \arrow[u, dotted] & \fbox{$x^2y^3$} \arrow[r, dotted] \arrow[u, dotted] & \cdots \\
x \arrow[u] \arrow[r] & xy \arrow[r] \arrow[u] & \fbox{$xy^2$} \arrow[r, dotted] \arrow[u, dotted] & \fbox{$xy^3$} \arrow[r, dotted] \arrow[u, dotted] & \cdots \\
1 \arrow[u] \arrow[r] & y \arrow[r] \arrow[u] & y^2 \arrow[r] \arrow[u] & y^3 \arrow[r] \arrow[u] & \cdots
\end{tikzcd}}
\medskip

For $\CC[X,Y]/(XY)$-modules this picture does not work quite so well since $XY = 0 = YX$. The more appropriate picture is a kind of "\emph{splitting}" of the lattice $\NN^2 \times \NN^2$, two copies since we are identifying with the space $\{0,1\}^{\omega} \amalg \{0,1\}^{\omega}$.

In Figure \ref{Sym4 and Sym6}, we see $\Sym^4(\CC^2)$ and $\Sym^6(\CC^2)$ represented as diagonals in the lattice $\NN^2$. These can be seen as the direct sums of the simple modules $\top(M) = M/\rad(M)$ of free modules $M \in \Mod(\CC[X,Y])$. In particular $\ttop(x^ny^m) \cong \CC$. We note, the spaces $\Sym^n(\CC^2)$ are invariant under the action of $\SL_2(\CC)$. Maps between free modules 
\[ f \in \Hom_{\CC[X,Y]}(M,N) \]
are always given by paths from $\ttop(N) \to \ttop(M)$ in the lattice $\NN^2$, modulo commutativity of the paths given by $XY=YX$. 

\medskip
\adjustbox{scale=0.7,center}{
\begin{tikzcd}\label{Sym4 and Sym6}
\fbox{$x^6$} \arrow[r] & x^6y \arrow[r] & x^6y^2 \arrow[r] & x^6y^3 \arrow[r] & x^6y^4 \arrow[r] & x^6y^5 \arrow[r] & x^6y^6 \\
x^5 \arrow[r] \arrow[u] & \fbox{$x^5y$} \arrow[r] \arrow[u] & x^5y^2 \arrow[r] \arrow[u] & x^5y^3 \arrow[r] \arrow[u] & x^5y^4 \arrow[r] \arrow[u] & x^5y^5 \arrow[r] \arrow[u] & x^5y^6 \arrow[u] \\
\fbox{$x^4$} \arrow[u] \arrow[r] \arrow[ru, no head, dotted] & x^4y \arrow[u] \arrow[r] & \fbox{$x^4y^2$} \arrow[u] \arrow[r] & x^4y^3 \arrow[u] \arrow[r] & x^4y^4 \arrow[u] \arrow[r] & x^4y^5 \arrow[u] \arrow[r] & x^4y^5 \arrow[u] \\
x^3 \arrow[u] \arrow[r] & \fbox{$x^3y$} \arrow[u] \arrow[r] \arrow[ru, no head, dotted] & x^3y^2 \arrow[r] \arrow[u] & \fbox{$x^3y^3$} \arrow[r] \arrow[u] & x^3y^4 \arrow[r] \arrow[u] & x^3y^5 \arrow[u] \arrow[r] & x^3y^6 \arrow[u] \\
x^2 \arrow[u] \arrow[r] & x^2y \arrow[u] \arrow[r] & \fbox{$x^2y^2$} \arrow[r] \arrow[u] \arrow[ru, no head, dotted] & x^2y^3 \arrow[r] \arrow[u] & \fbox{$x^2y^4$} \arrow[r] \arrow[u] & x^2y^5 \arrow[u] \arrow[r] & x^2y^6 \arrow[u] \\
x \arrow[u] \arrow[r] & xy \arrow[u] \arrow[r] & xy^2 \arrow[r] \arrow[u] & \fbox{$xy^3$} \arrow[r] \arrow[u] \arrow[ru, no head, dotted] & xy^4 \arrow[r] \arrow[u] & \fbox{$xy^5$} \arrow[u] \arrow[r] & xy^6 \arrow[u] \\
1 \arrow[u] \arrow[r] & y \arrow[r] \arrow[u] & y^2 \arrow[r] \arrow[u] & y^3 \arrow[r] \arrow[u] & \fbox{$y^4$} \arrow[r] \arrow[u] \arrow[ru, no head, dotted] & y^5 \arrow[r] \arrow[u] & \fbox{$y^6$} \arrow[u]
\end{tikzcd}}
\medskip

Taking the actions by $x$ and $y$ described above and consider them as maps
\[ \sym^n(\CC^2) \to \sym^{n+2}(\CC^2) \]
which can be thought of as the map given by the commutative square
\begin{center}
\begin{tikzcd}
 & x^{a+1}y^b \arrow[rd] &  \\
x^ay^b \arrow[ru] \arrow[rd] \arrow[rr, "(xy-yx)", no head, dashed] &  & x^{a+1}y^{b+1} \\
 & x^ay^{b+1} \arrow[ru] & 
\end{tikzcd}
\end{center}
 
We know that the ideals in $R = \CC^2[x,y]$ which are $\GL_2(\CC)$ fixed are powers of the maximal homogeneous ideal $\fm^n = (x,y)^n \cong \sym^n(\CC^2)$. We also know the action of 
 $\sln_2(\CC)$ on any $\sym^n(\CC^2)$ by taking the standard basis for $\sln_2(\CC)$, 
 \[ X = \begin{pmatrix}
 0 & 1 \\
 0 & 0
 \end{pmatrix}, \quad Y=\begin{pmatrix}
 0 & 0 \\
 1 & 0
 \end{pmatrix}, \quad H = \begin{pmatrix}
 1 & 0 \\
 0 & -1
 \end{pmatrix}\]
acting on the left by multiplication of vectors $v = \begin{pmatrix}
x \\ y
\end{pmatrix}$
 we see $X(v) = \begin{pmatrix}
 0 \\ y
 \end{pmatrix}$, $Y(v) = \begin{pmatrix}
 0 \\ x
 \end{pmatrix}$, and $H(v) = \begin{pmatrix}
 x \\ -y
 \end{pmatrix}$
This means, under the action of $\sln_2(\CC)$ on $\CC^2$ we get an action on $\CC[x,y]$ by defining 
\[ g \cdot f(x,y) = f(g(x,y)) \]

Remember, in \cite{GP}, Gel'fand and Ponomarev classified the indecomposable Harish-Chandra modules of the Lorentz group. They showed the classification of the finite dimensional modules was equivalent to classifying certain representation of $\sln_2(\CC)$ and of the Gel'fand-Ponomarev algebra $\CC[x,y]/(xy)$ given by the quiver with two loops as described at the beginning of this section. In particular, they showed that the indecomposable representations corresponded to "zig-zags" in the lattice $\NN^2$ with arrows pointing up and to the right as in the following diagram

\medskip

\begin{center}
\begin{tikzcd}
 \arrow[r] & \  &  &  &  \\
 &  \arrow[u] \arrow[r] \ & \  \arrow[r] & \ &  \\
 &  &  &  \arrow[u] \arrow[r] & \
\end{tikzcd}
\end{center}

\medskip

\subsection{Some Geometric Interpretations of The Commutative Algebra}

This should give us a clue as to what is actually happening in the case of surface algebras in general. How should we think of "reverse arrows"? In the ring of formal Laurent series
\[ \CC[[x,y]]/(xy) \subset \CC((x,y))/(xy) \]
the residues of $x$ and $y$ become invertible! One could think of such a thing as the structure sheaf of a fibre-product of two punctured disks 
\[\DD^{\times}_x \times \DD^{\times}_y = \Spec \CC((x)) \times \Spec \CC((y))\] 

given by the pullback
\[ \xymatrix{
\CC((x,y))/(xy) \ar[r] \ar[d] & \CC((x)) \times \CC((y)) \ar[d] \\
\CC = \left(\CC((x,y))/(xy)\right)/(x,y) \ar[r] & \CC \times \CC = \CC((x))/(x) \times \CC((y))/(y)
} \]
but there is more to the picture. In particular, in $\FF_x$ and $\FF_y$, $x$ and $y$ are of course invertible. 
Moreover, elements such as $x^py^q$ for $p, q \in \ZZ$ make sense in $\CC((x,y))/(xy)$, \emph{only when} $p \leq 0$ and $q \geq 0$, or vice-versa. Fascinating! \emph{This should point our noses in the right direction!} 

\begin{recall}
In \cite{S1}, the pullback of various "\emph{hereditary orders}" was defined. It was given by exactly such pullbacks on diagonal entries of various matrix algebras. 
\end{recall}

Suppose we take the algebraic completions of $\FF_x = \CC((x))$ and $\FF_y = \CC((y))$. What should happen? We obtain a new pullback,
\[\xymatrix{
\mathfrak{R} \ar[r] \ar[d] & \hat{\ZZ} \times \hat{\ZZ} \ar[d] \\
\CC \ar[r] & \CC \times \CC
}
\]
Here we should note, the Galois group of $\CC((x)) \times \CC((y))$ is
\[ \varprojlim \cG\left(\CC((x^{1/n}))/\CC((x))\right) \times \varprojlim \cG\left(\CC((y^{1/m}))/\CC((y))\right) \]
which is isomorphic to 
\[ \hat{\ZZ} \times \hat{\ZZ}.\]

But, \emph{what on earth} might $\mathfrak{R}$ be? Here is where one might want to think of deformation retracts
\[ \DD^{\times} \to S^1 \]
of punctured disks to circles. This is the meaning behind the gluings of the products of the matrix algebras in \cite{S1} given by the pullbacks of \emph{surface orders}, since 
$\GL_n(\FF)$ can be identified with maps $C^{\infty}(S^1, \GL_n(\CC))$. If we take the Lie algebra
\[ \fg((x)) \times \fg((y)) \]
of the two loop groups $G((x))$ and $G((y))$, we can take a Fourier transform,  
\[ \fg \otimes_{\CC} \CC((x)) \to \fg \otimes_{\CC} S^{1}\]
\[ \fg \otimes_{\CC} \CC((y)) \to \fg \otimes_{\CC} S^1\]
by 
\[ g_1 \otimes x^n \mapsto g_1 \otimes e^{-in\sigma_x} \]
and 
\[ g_2 \otimes y^m \mapsto g_2 \otimes e^{-im\sigma_y} \]
for $\sigma_x, \sigma_y \in [0,2\pi]$ giving coordinates on $S^1_x$ and $S^1_y$, the deformation retracts of $\DD^{\times}_x$ and $\DD^{\times}_y$ respectively. But then, reminding ourselves of the pullback, what does 
\[  \left(g_1 \otimes x^n\right) \times \left(g_2 \otimes y^m\right) \]
correspond to in the pullback algebra? What does
\[ \left(g_1 \otimes e^{-in\sigma_x}\right) \times \left(g_2 \otimes e^{-im\sigma_y}\right) \]
correspond to? 

\begin{recall}
For $\FF = \CC((x_j))$, we can associate an automorphism of $\FF((x_j^{1/n}))$, defined by 
\[ x_j \mapsto e^{2 \pi i k/n}\cdot x_j \]
to an element $\sigma_j \in \cG\left(\FF((x_j^{1/n_j}))/\CC((x_j))\right) \cong \ZZ/n_j\ZZ \hookrightarrow S_{N}$, where $n_j$ is the degree of the ramification at $x_j$. 
\end{recall}

But what should $k$ be? And why is this even important?

\begin{conj}
Whenever we take the pullback in this way, the least common multiple of the order of two such automorphisms should have a significant meaning in modular representations theory. In particular, if we are working over a coefficient ring or field $R$ of positive characteristic for the group ring $RG$, which divides the order of the group $G$, we get a decomposition of $RG = B_1 \oplus \cdots \oplus B_l$ into "blocks", which are indecomposable modules over $RG$. If our finite group $G$ is obtained from a pullback as above, and if the automorphisms $\sigma_r$ and $\sigma_s$ are mapped to primitive $r^{th}$ and $s^{th}$ roots of unity respectively, the least common multiple of the degrees of the ramifications over $x_r$ and $x_s$ yield a relation in the surface algebra, which is of the same type as that in the definition of Brauer graph algebras. This is how one obtains appropriate quotients of surface algebras yielding Brauer graph algebras.
\end{conj}

Again, to remind ourselves of why such a thing might be true, remember that the automorphisms as above correspond to a cyclic subgroups of the monodromy group of $X$, embedded as a permutation group in a symmetric group, and given as the pullback of the cyclic groups $\la \sigma_j \ra \cong \ZZ/n_j\ZZ$, over each ramification point/puncture $x_j$. Furthermore, the pullback of all such cyclic subgroups for $X$ is exactly the monodromy group given by the covering. This group can be seen as acting on the quiver of the surface algebra, as described in \cite{S1}, and it is was shown there that the two uniserial radicals of a projective module over a surface algebra "meet back up" at the least common multiple of the order of the corresponding permutations. The relations for Brauer graph algebras are exactly of this form, with some minor modifications allowed in addition.

\section{Loop Groups and Loop Algebras}

For most of the paper, we will make like simple and stick with the following simple conventional setup.

\subsection{The General Linear Group and its Lie Algebra}
$G = \GL_n(\CC)$ with be the \textbf{general linear group} of invertible $n \times n$ matrices over $\CC$, we will take the two \textbf{Borel subgroups}

\[ B_+ = \begin{pmatrix}
b_{1,1} & b_{1,2} & \cdots & b_{1,n} \\
0 & b_{2,2} & \cdots & b_{2,n} \\
\vdots & \vdots & \ddots & \vdots \\
0 & 0 & \cdots & b_{n,n}
\end{pmatrix}, \quad B_{-} =  \begin{pmatrix}
b_{1,1} & 0 & \cdots & 0 \\
b_{2,1} & b_{2,2} & \cdots & 0 \\
\vdots & \vdots & \ddots & \vdots \\
b_{n,1} & b_{n,2} \cdots & b_{n,n}
\end{pmatrix} \]

We will take $U_{+} = [B_+, B_+]$ to be the \textbf{unipotent radical} of $B_+$ (i.e. the subgroup of upper triangular matirces with $1$s on the diagonal) and we take $U_{-} \subset B_{-}$ to be the \emph{lower triangular} matrices with $1$s on the diagonal. We will take $T$ to be the standard torus, i.e. the subgroup with entries only on the diagonal. 

We will take $\fg = \gl_n$ to be the Lie algebra of $G$, $\fb_{\pm}$ the Lie algebras of $B_{\pm}$, $\fn_{\pm}$ will be the (nilpotent) Lie algebras of $U_{\pm}$, and $\fh$ will be the Lie algebra of $T$ (often called the \textbf{Cartan subalgebra} with respect to $T$). 

We can phrase many results in general terms for complex reductive algebraic groups, but would like to be able to give some concrete examples as we go that will give some visual representation of what is happening, and this can sometimes be more difficult if we work in full generality.

\subsection{Loop Groups}

\section{Complex Curves and the Heart of the Answer to the Correspondence}

Let $X$ be a smooth curve over $\CC$, and let $x \in X$. Let $\cO(X)$ be the structure sheaf, or if we wish to stay down to earth and if we are able to think of things in terms of affine algebraic varieties, we can take $\cO(X) = \CC[X]$ to be the coordinate ring of $X$. For an open set $U \in X$ we take $\cO(U)$ to be the sheaf of regular functions on $U$. Let $\cO(X)_{x}$ be the sheaf local rings with respect to $x \in X$, with corresponding maximal ideal $\fm_x$, or if it feels more concrete, think $\CC[X]_x$ to be the coordinate ring localized with respect to $x \in X$. We will use the notation $\O(X)_x$ to mean the sheaf of complete local rings, which one could think of locally as a copy of $\CC[[x]]$, i.e. power series in $x$. Let $\cK_x$ be the sheaf of the fields of fractions, which locally looks like $\CC((t))$, the Laurent series rings. We then have
\[ \Spec \pO_x \cong \DD_x \]
is a copy of the open disk, and 
\[ \Spec \cK_x \cong \DD_x^{\times} \]
is a copy of the punctured disk.

\section{Surface Algebras}
We can now introduce the \textit{surface algebras}, which along with their $\fm$-adic completions will be the main objects of study in what follows ($\fm$ being the arrow ideal).

Let $Q = (Q_0, Q_1, h, t)$ be a quiver, with the set of vertices $Q_0$, and the set of arrows $Q_1$. There are two maps, 
\[ t, h: Q_1 \to Q_0 \]
taking an arrow $a \in Q_1$ to its \textbf{head} $ha$, and \textbf{tail} $ta$. This is a refinement of the incidence map for an undirected graph, and we define 
\[ \partial a = \{\partial_{\bullet}a, \partial^{\bullet}a\}:= \{ta, ha\}.\] 
In this case the order is not arbitrary as it would be for undirected graphs. The \textbf{path algebra} of a quiver $Q$, denoted $kQ$, over a field $k$, is the $k$-vector space spanned by all oriented paths in $Q$. It is an associative algebra, and is finite dimensional as a $k$-vector space if and only if $Q$ has no oriented cycles. There are trivial paths $i \in Q_0$, given by the vertices, and multiplication in the path algebra is defined by concatenation of paths, when such a concatenation exists. Otherwise the multiplication is defined to be zero. More precisely, if $p$ and $q$ are directed paths in $Q$, and $hp=tq$, then $qp$ is defined as the concatenation of $p$ and $q$. Note, we will read paths from \emph{right to left}. Let $A=kQ$. The \textbf{vertex span} $A_0=k^{Q_0}$, and the \textbf{arrow span} $A_1=k^{Q_1}$ are finite dimensional subspaces. $A_0$ is a finite dimensional commutative $k$-algebra, and $A_1$ is an $A_0$-bimodule. The path algebra then has a grading by path length, 
\[ A = A_0 \la A_1 \ra = \bigoplus_{d=0}^{\infty} A^{\otimes d}.\]
The path algebra $A$ has primitive orthogonal idempotents $\{e_i\}_{i \in Q_0}$. Let $A_{i,j} = e_jAe_i$ be the $k$-linear span of paths in $Q$, from vertex $i$ to $j$. Let $\fm = \prod_{d=1}^{\infty} A^{\otimes d}$ denote the \emph{arrow ideal} of $Q$, generated by the arrows $Q_1$. We will define the \textbf{complete path algebra} to be
\[ \cA = A_0 \la \la A_1 \ra \ra = \prod_{d=0}^{\infty} A^{\otimes d} .\]
We put the $\fm$-adic topology on $\cA$, with neighborhoods of $0$ generated by $\fm^n$. The elements of $\cA$ are all formal linear combinations of paths, including infinite linear combinations. If $\phi: \cA \to \cA$ is an automorphism fixing $\cA_0$ then $\phi$ is continuous in the $\fm$-adic topology, and $\fm$ is invariant under such algebra automorphisms.

\begin{defn}
An \textbf{ideal} in the path algebra $A$ will be a two sided ideal generated by linear combinations of paths which share a common starting vertex and terminal vertex in the quiver. The \textbf{quotient path algebra} of a quiver with relations will be the quotient by this ideal. 
\end{defn}

Next let us turn to the specific quivers with relations of interest for our current purposes. 

\begin{defn}
Define a \textbf{free surface algebra} to be the path algebra of the medial quiver of any combinatorial map (i.e. a cellularly embedded graph) given by a constellation $C=[\sigma, \alpha, \phi]$.
\end{defn}

\begin{defn}\label{surface algebra}
	Let $Q = (Q_0, Q_1)$ be a finite connected quiver. Then we say the bound path algebra $\Lambda = kQ/I$ is a \textbf{surface algebra} if the following properties hold:
	
	\begin{enumerate}
		\item For every vertex $x \in Q_0$ there are exactly two arrows $a, a' \in Q_1$ with $ha = x = ha'$, and exactly two arrows $b, b \in Q_1$ such that $tb = x = tb'$.
		\item For any arrow $a \in Q_1$ there is exactly one arrow $b \in Q_1$ such that $ba \in I$, and there is exactly one arrow $c \in Q_1$ such that $ac \in I$.
		\item For any arrow $a \in Q_1$ there is exactly one arrow $b' \in Q_1$ such that $b'a \notin I$, and there is exactly one arrow $c' \in Q_1$ such that $ac' \notin I$.
		\item The ideal $I$ is generated by paths of length $2$. 
	\end{enumerate}
\end{defn}

\section{Surface Orders}
In this section, take $\FF$ to be the field of fractions of the commutative ring $R$, which is a complete discrete valuation ring. Let $\mathfrak{m}$ denote the unique maximal ideal, and $k = R/\mathfrak{m}$ the residue field. Concrete examples of the setup which are extremely important are
\[
\begin{matrix}
\FF = \CC((x)) & R = \CC[[x]] & \fm = (x) & \CC = R/\fm \\
\FF = \hat{\QQ}_{(p)} & R = \hat{\ZZ}_{(p)} & \fm = (p) & k = \ZZ/(p) = \FF_p \\
\FF = \FF_p((x)) & R = \FF_p[[x]] & \fm = (x) & k = \FF_p
\end{matrix} \]

We stay mostly with the case of $\FF = \CC((x))$, the field of formal Laurent series, but for now let us define surface orders in the generality of \cite{S1}. 

\begin{defn}
An $R$-\textbf{lattice} is a finitely generated projective module over $R$. In particular, if $R$ is a Dedekind domain, every $R$-lattice is finitely generated and torsion free. 
\end{defn}

\begin{ex}
For example, if $R = \ZZ$, then $\ZZ^2$ is a $\ZZ$-lattice via addition of ordered pairs. As another example, let $R = \CC[x, y]/(xy)$. Then $\CC[x.y]/(xy)$ is an $R$-lattice over itself via the action given by multiplication by $\overline{x}$ and $\overline{y}$, the residues of $x$ and $y$ in $R$. One can visualize this via the maps of the bigraded shifts
\[
\xymatrix{
& R \ar[dl]_{x} \ar[dr]^y & \\
R(-1,0) \ar[d]_x & & R(0,-1) \ar[d]^y \\
R(-2,0) \ar[d]_x & & R(0,-2) \ar[d]^y \\
R(-3,0) \ar[d]_x & & R(0,-3) \ar[d]^y \\
\vdots & & \vdots
}
\]

\end{ex}
	
\begin{defn}
An $R$-\textbf{Order} $\Lambda$ in a $k$-algebra $A$ is a unital subring of $A$ such that
	\begin{enumerate}
		\item $\FF \Lambda = A$, and
		\item $\Lambda$ is finitely generated as an $R$-module. 
	\end{enumerate}
\end{defn}

\begin{defn}
	Let $C = [\sigma, \alpha, \phi]$ be a \textbf{constellation}\footnote{Recall from \cite{S1}, a constellation is a permutation group acting transitively on $[2N]$, with $\alpha$ a fixed point free involution, and $\sigma \alpha \phi =1$.}, and let $\Gamma \hookrightarrow \Sigma$ be the associated graph cellularly embedded in the compact Riemann surface $\Sigma$. Further, let $n(i)=n_i = |\sigma_i|$ denote the length of the cycle $\sigma_i$ in the permutation $\sigma = \sigma_1 \sigma_2 \cdots \sigma_p$. Remember, for a constellation $C$, and the associated graph $\Gamma$, the length of the (nonzero) cycle in the medial quiver $Q(C)$ with gentle relations, which is associated to $\sigma_i$ (and its corresponding vertex of $\Gamma$) is just the order of the cycle $\sigma_i$.
	\begin{enumerate}
		\item For each cycle $\sigma_i$, corresponding to the vertex $i \in \Gamma_0$, we associate a Dedekind domain $R_i$, with a maximal ideal $\fm_i$, and a \textbf{vertex} ($R_i$-)\textbf{order}, $\Lambda(\sigma_i) = \Lambda_i$. 
		\item The \emph{vertex order} associated to the cycle $\sigma_i$ is then given by a matrix $R$-subalgebra of $\Mat_{n_i \times n_i}(R_i)$. 
		\[
		\Lambda_i = \begin{pmatrix}
			R_i & \fm_i & \fm_i & \cdots & \fm_i & \fm_i \\
			R_i & R_i   & \fm_i & \cdots & \fm_i & \fm_i \\
			R_i & R_i   & R_i   & \cdots  & \fm_i & \fm_i \\
	   \vdots & \vdots & \vdots & \ddots & \vdots & \vdots \\
	   R_i & R_i & R_i & \cdots & R_i & \fm_i \\
	   R_i & R_i & R_i & \cdots & R_i & R_i 
		\end{pmatrix}_{n(i)}
		\]
		If we take $R_i$ to be $k[[x_i]]$, with maximal ideal $(x)$, we have that $\Lambda_i$ is a \emph{hereditary order} \footnote{Recall, an algebra, or an order, is \textbf{hereditary} if no module has a minimal projective resolution greater than length one. This means the \textit{projective dimension} of any module is no greater than one, and therefore the global dimension of the algebra is at most one.}.		
		
		\item Let $\Lambda_i^{(k,k)}$ denote the $(k,k)$ entry of $\Lambda_i$ (in $R_i$).  
		\item For each $1 \leq k \leq n_i$ let
		\[ 
		P_{i,1}:= \begin{pmatrix}
		R_i \\ \vdots \\ R_i \\ R_i \\ R_i \\ \vdots \\ R_i \\ R_i
		\end{pmatrix},		 
		\ P(\sigma_i) = \sigma_i \cdot P_{i,1}:= \begin{pmatrix}
		\fm_i \\ R_i \\ \vdots \\ R_i \\ R_i \\ \vdots \\ R_i \\ R_i
		\end{pmatrix}, \cdots,
		\ P(\sigma_i^{k-1}) = P_{i,k}= \begin{pmatrix}
		\fm_i \\ \vdots \\ \fm_i \\ R_i \\ R_i \\ \vdots \\ R_i \\ R_i
		\end{pmatrix}, \cdots ,
		P_{i,n_i}:= \begin{pmatrix}
		\fm_i \\ \vdots \\ \fm_i \\ \fm_i \\ \fm_i \\ \vdots \\ \fm_i \\ R_i
		\end{pmatrix} 
		\]
		where the $k^{th}$ entry is the first entry equal to $R_i$ for $P_{i,k} = \sigma_i^{k-1} \cdot P_{i,1} = P(\sigma_i^{k-1})$. 
	\end{enumerate}
\end{defn}

The modules $\{P_{i,k}: \ 1 \leq k \leq n_i\}$ give a complete set of non-isomorphic indecomposable projective (left) $\Lambda_i$-modules, with the natural inclusions
\[ P_{i,1} \hookleftarrow P_{i,2} \hookleftarrow \cdots \hookleftarrow P_{i,n_i-1} \hookleftarrow P_{i,n_i} 
\hookleftarrow P_{i,1}. \]
where the final map is given by left-multiplication by $\fm_i$. If we identify $P_{i,k}$ with the edge 
$e_k^i = e_k(\sigma_i)$, where $\sigma_i = (e_1^i, e_2^i, ..., e_{n_i}^i)$ is a cyclic permutation, then the chain of 
inclusions can be interpreted in terms of the cycle $\sigma_i$. From the embedding $\Gamma
\hookrightarrow \Sigma$ given by the constellation $C = [\sigma, \alpha, \phi]$, this can be interpreted as 
walking clockwise around the vertex of $\sigma_i$. We will take $P_{i,k} = P_{i,k+n_i}$, and each $e_k^i$ is 
multiplied by the automorphism $\sigma_i$, i.e. there is some multiplication by a 
power of $\fm_i$ involved. 
\indent In particular, multiplication by 
\[ \sigma_i= \begin{pmatrix}
	0 	  & 0 & 0 & \cdots & 0 & \fm_i \\
	1 	  & 0 & 0 & \cdots & 0 & 0 \\
	0 	  & 1 & 0 & \cdots & 0 & 0 \\
\vdots & \vdots & \vdots & \ddots & \vdots & \vdots \\
	0 	  & 0 & 0 &\cdots & 0 & 0 \\
	0 	  & 0 & 0 &\cdots & 1 & 0 
\end{pmatrix}_{n_i} \]
cyclically permutes the indecomposable projective $\Lambda_i$-modules $P_{i,k}$, and it induces an automorphism of the matrix algebra $\Lambda_i$ which we also call $\sigma_i$. Now, for each pair of cycles $\sigma_i, \sigma_j \in S_{[2m]}$ of $\sigma$, we fix an 
isomorphism
\[ R_i/\fm_i \cong R_j/\fm_j. \] 
Identifying all such rings, let $k = R_i/\fm_i$ for all $\sigma_i \in \Gamma_0$. Let $
\pi_i: R_i \to k$ be a fixed epimorphism with kernel $\fm_i$ a maximal ideal of $R_i$. Now, we have a pull-back 
diagram
\[ 
\xymatrix{
R_{i,j} \ar[r]^{\tilde{\pi}_i}  \ar[d]_{\tilde{\pi}_j} & R_i \ar[d]^{\pi_i} \\
R_j \ar[r]_{\pi_j} & k 
}
\]
which is in general different and non-isomorphic for different choices of $\pi_i$ and $\pi_j$.

\begin{defn}
	Let $\cN(\Lambda) = \prod_{\sigma_i \in \Gamma_0} \Lambda_i$. Let $e_k^i$ be an edge around $\sigma_i$, and let $
	\alpha^{i,j}_{k,l} = (e_k^i, e_l^j)$ be a $2$-cycle of the fixed-point free involution $\alpha$ of $C=[\sigma, \alpha, \phi]$ giving the end vertices $
	\sigma_i$ and $\sigma_j$ of the edge $e_k^i \equiv e_l^j$ under the gluing identifying the half-edges $e_k^i$ and 
	$e_l^j$. It is possible that $\sigma_i=\sigma_j$ if $\alpha^{i,j}_{k,l}$ defines a loop at the vertex $\sigma_i$ in $
	\Gamma$. We replace the product $\Lambda_i^{(k,k)} \times \Lambda_j^{(l,l)}$ in $\Lambda_i \times \Lambda_j$ with $R_{i,j}$. 
	This identifies the $(k,k)$ entry of $\Lambda_i$ with the $(l,l)$ entry of $\Lambda_j$, modulo $\fm_{i,j}$, the ideal of $R_{i,j}$ given via the pullback of $\fm_i$ and $\fm_j$. Doing this for all 
	edges of $\Gamma$, we get the \textbf{Dessin Order} $\Lambda := 
	\Lambda(C) = \Lambda(\Gamma)$ associated to the constellation $C$, or equivalently to the embedded 
	graph $\Gamma \hookrightarrow \Sigma$. We will call the hereditary order $\prod_{\sigma_i}\Lambda_i$ the (noncommutative)\textbf{normalization} of the dessin order order $\Lambda$. 
\end{defn}

		The indecomposable projective $\Lambda$-modules are in bijection with the $2$-cycles $(e^i_k, e^j_l)=\alpha^{i,j}_{k,l}$, of $\alpha \in S_{2m}$ for the constellation $C=[\sigma, \alpha, \phi]$. Equivalently, the indecomposable projectives are in bijection with the edges $\Gamma_1$. We label them as $P_{e}$ for $\alpha^{i,j}_{k,l} = e = (e^i_k, e^j_l) \in \Gamma_1$ attached to the vertices $\sigma_i$ and $\sigma_j$.

\section{Affine Grassmannians and Affine Flag Varieties}

In this section, we take follow Magyar \cite{M} and \cite{F} to give the setup and definitions of affine Grassmannians and affine Flag Varieties. Let $\FF = \CC((x))$ be the field of formal Laurent series, and $R = \CC[[x]]$ the completion of the polynomial ring $\CC[x]$ (i.e. formal Taylor series) with respect to $\fm = (x)$. Since we will work with \emph{formal} power series and \emph{formal} Laurent series, we are allowing infinite sums, and evaluation at $x=0$ in general is the only evaluation allowed. Moreover, we define the \textbf{formal disk} to be
\[ D^{\times} = \Spec \CC((x)) \]

Let $\GL_n(\FF) = G((x))$, be the group of invertible matrices over $\CC((x))$, i.e. the \textbf{loop group}. Define a topology on $\GL_n(\FF)$ by letting the base open neighborhoods of $\id \in \GL_n(\FF)$ to be given by congruence subgroups $K_N$, for $N \in \ZZ_{\geq 0}$. In particular, we take 
 
\[ G_j := \{g \in G: \ \ord \det(g) = j\} \]
where $\ord(f)$ for $f = \sum_{i>N}a_ix^i \in \CC((x))$ is the smallest $N$ such that $a_N \neq 0$. In particular, the elements in $\FF$ are only allowed to have finitely many negative terms, so that we do not get "\emph{poles}" of infinite order. Then we have the decomposition 
\[ \GL_n(\FF) = \amalg_{j \in \ZZ} G_j.\]

Taking $\sigma \in G_1$ we have $\sigma^j \cdot G_0 = G_0 \cdot \sigma^j = G_j$. 

\begin{defn}
Next, let $V$ be a vector space over $\FF = \CC((x))$. If $V = \FF^n$ then of course $V$ is finite dimensional over $\FF$. The group $\GL_n(\FF)$ acts on $V$. In particular, taking $\{e_1, ..., e_n\}$ to be the standard basis (over $\FF$), for any $r \in \ZZ$ define
\[ e_{i+rn} = x^r \cdot e_i.\]
This gives a $\CC$-basis $\{v_i\}_{i \in \ZZ} = \{x^r \cdot e_i\}_{r \in \ZZ}^{i=1,...,n}$ of $V$. So we may take 
\[ \Lambda = Rv_1 \oplus Rv_2 \oplus \cdots \oplus Rv_n \]
to be the $R$-submodule in $V$ spanned by some $\FF$-basis $\{v_1, ..., v_n\}$ of $V$ from $\{v_i\}_{i \in \ZZ}$. This gives
\[ \Lambda = \Span_{\CC} \la x^{r_i} \cdot e_i \ra_{r \geq 0}^{i=1,...,n} \]
and this gives $\Lambda$ an $R$-\emph{lattice structure}. 

\begin{remark}
Here we are allowing infinite linear combinations over $\CC$ of vectors $e_{i+rn} = x^r \cdot e_i$, but with coefficients in $\FF$, $\Lambda$ is a \emph{finite $\FF$-dimensional vector space}. It is obviously a \emph{finitely generated projective module over $R$} and it is therefore an $R$-lattice. 
\end{remark}

\begin{recall}
An $R$-\textbf{Order} $\Lambda$ in a $\CC$-algebra $A$ is a unital subring of $A$ such that
	\begin{enumerate}
		\item $\FF \Lambda = A$, and
		\item $\Lambda$ is finitely generated as an $R$-module. 
	\end{enumerate}
For example, 
\[ \Lambda = \begin{pmatrix}
R & \fm \\
R & R
\end{pmatrix} \subset 
\Mat_{2 \times 2}(\FF) = \FF \otimes_{\CC} \Mat_{2 \times 2}(\CC) \subset \Mat_{2 \times 2}(R) = R \otimes_{\CC} \Mat_{2 \times 2}(\CC). \]
is an $R$-order in the $\CC$-algebra $A = \Mat_{2 \times 2}(\FF)$. 
\end{recall}

Now, let $\{e_i\}_{i=1}^n$ be the standard $\FF$-basis. Define $\sigma \cdot e_i = e_{i+1}$ to be the \textbf{shift operator} given by the matrix
\[ \sigma = \begin{pmatrix}
0 & 0 & 0 & \cdots & 0 & x \\
1 & 0 & 0 & \cdots & 0 & 0 \\
0 & 1 & 0 & \cdots & 0 & 0 \\
\vdots & \vdots & \vdots & \ddots & \vdots & \vdots \\
0 & 0 & 0 & \cdots & 1 & 0
\end{pmatrix} \]
\end{defn}

\begin{defn}
We define the \textbf{affine Grassmannian} $\Gr(V)$ to be the space of all $R$-lattices in $V$. To be formal, we can let $M$ be an $R$-module with an ordered set $m_1, ..., m_k \in M$. Define an equivalence relation
\[ (M; m_1, ..., m_k) \sim (M'; m_1', ..., m_k') \]
if and only if there is an isomorphism $f: M \to M'$ such that $f(m_i) = m_i'$. Define the \textbf{Grassmanian functor}, denoted $\Gr(r, k)$, to be the functor from the category of rings to the category of sets which assigns to a $k$-algebra $R$
\[ \left\lbrace \substack{ (M;m_1, ..., m_k), \ \rank(M)=r, \\ M=R\la m_1, ..., m_k\ra } \right\rbrace \bigg/ \sim\]
and to each ring homomorphism $f:R \to S$ it assigns a set map
\[ [M;, m_1, ..., m_k] \mapsto [M \otimes_R S; m_1 \otimes 1, ..., m_k \otimes 1],\]
where $[M; m_1, ..., m_k]$ denotes the equivalence classes of locally free $R$-modules. We will call the lattice 
\[ E_{\bullet} = Re_1 \oplus Re_2 \oplus \cdots Re_n = E_1 \oplus E_2 \oplus \cdots \oplus E_n \]
the \textbf{standard lattice}. The action of $\sigma_x$ gives $E_i \mapsto E_{i+1}$ and so $\sigma_x^r \cdot E_i = E_{i+rn}$.

Taking
\[ E_{\bullet}(0) = \Span_R \la e_1, ..., e_n \ra \]
to be the standard lattice, then a $\CC$-basis for $E_{\bullet}(0)$ is 
\[ \Span_{\CC} \{\sigma_x^i e_1, \sigma^i_x e_2, ...., \sigma_x^i e_n \}_{i \geq 0}\]
 as before, and define
\[ E_{\bullet}(i,j) = \Span_R \{\sigma^i_x \cdot e_j, \sigma^i_x \cdot e_{j+1}, ..., \sigma^i_x \cdot e_{j+n-1}\}\]
One might think of this as first some \emph{"rotation"} of the basis, where \emph{"wrapping around"} back to 
$e_1$ induces multiplication by $x$, followed by an application of $\sigma_x^i$ which \emph{"shifts"} every $\CC
$-degree by $i$. In particular, if we have the $R$-lattice $E_{\bullet}(i,j)$ its basis over $\CC$ is 
$\Span_{\CC} \{e_i\}_{i \geq j}$ taking the convention $\sigma_x^i \cdot E_j = E_{j+in} = \Span_R x^ie_j$.

We would like to understand the stabilizer of $E_{\bullet}(i,j)$. 

If $P_0 := \GL_n(R) \subseteq \GL_n(\FF)$, is the subgroup of matrices over $R = \CC[[x]]$ with $\ord(\det(g)) = 0$, then we have $\Gr(V) \cong \GL_n(\FF)/P_0$ as $P_0$ stabilizes
\[ E_{\bullet}(0) = \Span_R\{e_1, ..., e_n\} = \Span_{\CC}\{x^r \cdot e_1, ..., x^r \cdot e_n\}_{r \geq 0} \]

The connected components of $\Gr(V)$ are
\[ \Gr_j(V) = G_0 \cdot E_{\bullet}(j) = G_j \cdot E_{\bullet}(0) \cong G_j/P_0 \]

The \textbf{complete affine flag variety}, $\Fl(V)$ is the space of all flags of $R$-lattices
\[ \Lambda_{\bullet} = \Lambda_1 \subset \Lambda_2 \subset \cdots \subset \Lambda_n \]
with $x \cdot \Lambda_i \subset \Lambda_{i+1}$.  
\end{defn}

Now, we have the following diagram
\[
\xymatrix{
\FF = \CC((x)) & \\
R = \CC[[x]] \ar[u]_{\subseteq} \ar[r]_{\pi}^{f(x) \mapsto f(0)} & \CC = \CC[[x]]/(x)
}\]
which gives the diagram
\[
\xymatrix{
G = \GL_n(\FF) & \\
P_0 = \GL_n(R) \ar[u]_{\subset} \ar[r]^{\pi}_{f(x) \mapsto f(0)} & \GL_n(\CC) \\
B_+ = \pi^{-1}(B(\CC)) \ar[u]_{\subset} \ar[r]^{\pi}_{f(x) \mapsto f(0)} & B(\CC)_+ \ar[u]_{\subset} 
}\]
where $B(\CC)_+$ is a \textbf{Borel subgroup} of \emph{upper triangular matrices} in $\GL_n(\CC)$. The subgroup $B_+$ is generally called the \textbf{Iwahara subgroup} of $G = \GL_n(\FF)$. We then have that the \textbf{flag variety} can be realized as
\[ \GL_n(\CC)/B(\CC)_+ \]
and the \textbf{affine flag variety} can be realized as
\[ G/B = \GL_n(\FF)/B_+ \]
and then $G/P_0 = \GL_n(\FF)/\GL_n(R)$ is as before, the affine Grassmannian (sometimes also called the \textbf{loop Grassmannian}). We have that $\dim_{\FF}(\Lambda_j/\Lambda_{j+1}) = 1$ for lattices in the complete affine flag variety.

Now, we have that the affine flag variety
\[ \Fl(V) \cong G/B = \GL_n(\FF)/B_+ \]
has connected components 
\[ \Fl_j(V) = G_j \cdot E_{\bullet}(0) = G_0 \cdot E_{\bullet}(j) \cong G_j/B_+\]
In particular, if we define
\[ \Lambda_{\bullet}(0) = E_1 \subset \left(E_1 \oplus E_2\right) \subset \cdots \subset \left(\bigoplus_{i=1}^n E_i \right) \]
to be the \textbf{standard flag}, we can take 
\[ \Lambda_{\bullet}(i,j) = x^i \cdot E_j \subset x^i\left(E_j \oplus E_{j+1})\right) \subset \cdots \subset x^i\left( \bigoplus_{k=j}^{j+n-1} E_k \right)\]
The stabilizer of $\Lambda_{\bullet}(i,j)$ is the subgroup of matrices
\[B(j)_+ = \{A \in \GL_n(R): \ \deg(a_{p,q}) > j \ \forall p<q \]
So for example $B(r)_+$ would be
\[ \begin{pmatrix}
x^r \CC[[x]] & x^r \CC[[x]] & \cdots & x^r \CC[[x]] \\
\CC[[x]] & \CC[[x]] & \cdots & x^r \CC[[x]] \\
\vdots & \vdots & \ddots & \vdots \\
\CC[[x]] & \CC[[x]] & \cdots & \CC[[x]]
\end{pmatrix} \]
gives lower triangular matrices over $R/(x^r)$. Moreover, if we let $(x) = \fm$, in our earlier notation $B(1)_+$ is 
\[ \begin{pmatrix}
			R_i & \fm_i & \fm_i & \cdots & \fm_i & \fm_i \\
			R_i & R_i   & \fm_i & \cdots & \fm_i & \fm_i \\
			R_i & R_i   & R_i   & \cdots  & \fm_i & \fm_i \\
	   \vdots & \vdots & \vdots & \ddots & \vdots & \vdots \\
	   R_i & R_i & R_i & \cdots & R_i & \fm_i \\
	   R_i & R_i & R_i & \cdots & R_i & R_i 
		\end{pmatrix}_{n(i)}\]
and if we take its Lie algebra it is of the same form. This Lie algebra can be identified with $\Lambda_{\bullet}(1)$. We have $\Fl(V) \cong \GL_n(R)/B$.

\section{The Affine Weyl Groups}
Let $G((x)) = \GL_n(\FF)$ be the loop group with $\FF = \CC((x))$, and let $\GL_n(\CC)$ be the usual group of invertible matrices over $\CC$. Let $\fg = \mathfrak{gl}_n(\CC)$ be its Lie algebra. Take the \textbf{loop algebra}
\[ \mathfrak{g}((x)) = \mathfrak{gl}_n(\CC) \otimes_{\CC} \CC((x)) \]
As in the description of the Gel'fand-Ponomarev Algebra in Section \ref{Gel'fand Ponomarev}, we may take the Fourier transform of $\fg((x))$ to get the group
\[ \gl_n(\CC) \otimes_{\CC} S^1 \]
via 
\[ g \otimes_{\CC} x^n \mapsto g \otimes_{\CC} e^{-in\sigma_x} \]
for $\sigma_x \in [0, 2\pi]$. 

The Weyl group of $\GL_n(\CC)$ is of course the symmetric group $W = S_n$. We would like to define the Weyl group for $G((x))$ now. First, let $S_{\infty}$ ($n \geq 2$) to be all bijections $s: \ZZ \to \ZZ$. Let $\sigma(i) = i+1$ be the shift operator on $\ZZ$. So, $\sigma^n(i) = i+n$. As in \cite{BB}, define the subgroup $\tilde{S}_n \subset S_{\infty}$ to be the subgroup such that 

\begin{enumerate}
\item elements $\tilde{s} \in \tilde{S}_n$ satisfy $\tilde{s}(r+n) = \tilde{s}(r)+n$ for all $r \in \ZZ$. In other words $\tilde{s}$ commutes with $\sigma^n$. 
\item $\sum_{r=1}^n \tilde{s}(r) = {n+1 \choose 2}$. 
\end{enumerate}
The elements $\tilde{s}$ of the group $\tilde{S}_n$ are completely determined by their values on $[n]$ in the following way. Let
\[ \begin{pmatrix}
1 & 2 & \cdots & n \\
\tilde{s}(1) & \tilde{s}(2) & \cdots & \tilde{s}(n)
\end{pmatrix} = [a_1, a_2, ..., a_n] \]
be represented by the $n$-tuple $[a_1, a_2, ..., a_n] \in \ZZ^n$, where $a_i = \tilde{s}(i) \in \ZZ$ for $i=1, ..., n$. We follow the terminology of \cite{BB} and call this the \textbf{window} for $s$, or $s$ in the \textbf{window notation}. As generators we may take
\[ \tilde{s}_i = [1, 2, ..., i-1, i+1, i, i+2, ..., n] \]
which can be uniquely identified with the usual generator $s_i = (i, i+1)$ of $S_n$, given by the simple reflection and corresponding to the $i^{th}$-vertex of the $A_n$ Coxeter graph
\[ \xymatrix{
\bullet_1 \ar@{-}[r] & \bullet_2 \ar@{-}[r] & \cdots & \bullet_{i-1} \ar@{-}[r] & \bullet_{i} \ar@{-}[r] & \bullet_{i+1} \ar@{-}[r] & \cdots & \bullet_{n-2} \ar@{-}[r] & \bullet_{n-1}
}\]
We may think of $\tilde{s}_i \in \tilde{S}_n$ as the "\emph{affine version}" of $S_n$, given by looking through the "$\ZZ/n\ZZ \times \ZZ/n\ZZ$ windows" of $\ZZ \times \ZZ$. The only difference that would be prudent to point out is that this means
\[ \tilde{s}_n = [0, 2, 3, ..., n-1, n+1]. \]
In Magyar's setup, this is explained as follows. Take some $a = [a_1, ..., a_n] \in \ZZ^n$, and define
\[ \tau^a(i) = \tau^{[a_1, ..., a_n]}(i) = i+na = i+n[a_1, ...., a_n] \in \ZZ, \]
where we take $i \in \ZZ/n\ZZ$. Then for all $\tilde{s} \in \tilde{S}_n$ we can write $\tilde{s} = s\tau^{[a_1, ..., a_n]}$, where $s \in S_n$. This gives a semi-direct product structure on $\tilde{S}_n$ via the embeddings
\[ \xymatrix{
\ZZ^n \ar[r]^{\triangleleft} & \tilde{S_n} \\
S_n \ar[ur]_{\subset} & 
}\]
and $\tilde{S}_n = S_n \rtimes \ZZ^n$. Let $V = \FF^n$ have $\FF$-basis $\{e_1, ..., e_n\}$, and let $\tilde{s} \in \tilde{S}_n$ act by $\tilde{s} \cdot e_j = e_{\tilde{s}(j)}$. Now, as we have already mentioned, we may think of the action of $\tilde{S}_n$ on $\ZZ$ as restricting to an action on $\ZZ/n\ZZ$. But, what is more elucidating (and fascinating) is that we can view this as an action on the points $e^{2 \pi i k/n} \in S^1$ via the \emph{wiring diagrams} of Bernstein, Fomin, and Zelevinsky \cite{BFZ}. Their construction is to think of the points
\[ (p, e^{2 \pi i k/n}) \in [0,1] \times S^1 \]
where for $p = 0$, we can label the discrete subset $(0, e^{2 \pi i k/n}) = (0,k) \in [0,1] \times \ZZ/n\ZZ$. Next, choosing some $a = [a_1, a_2, ...., a_n] \in \ZZ^n$, we map $k \mapsto s(k)$, but the path $\gamma_s \subset [0,1] \times S^1$, which connects $k \to s(k)$, wraps around the cylinder $a_k$ times. The picture we should in mind is Figure \ref{Affine Permutation}:

\medskip

\begin{center}\label{Affine Permutation}
\includegraphics[
    page=1,
    width=200pt,
    height=200pt,
    keepaspectratio]
    {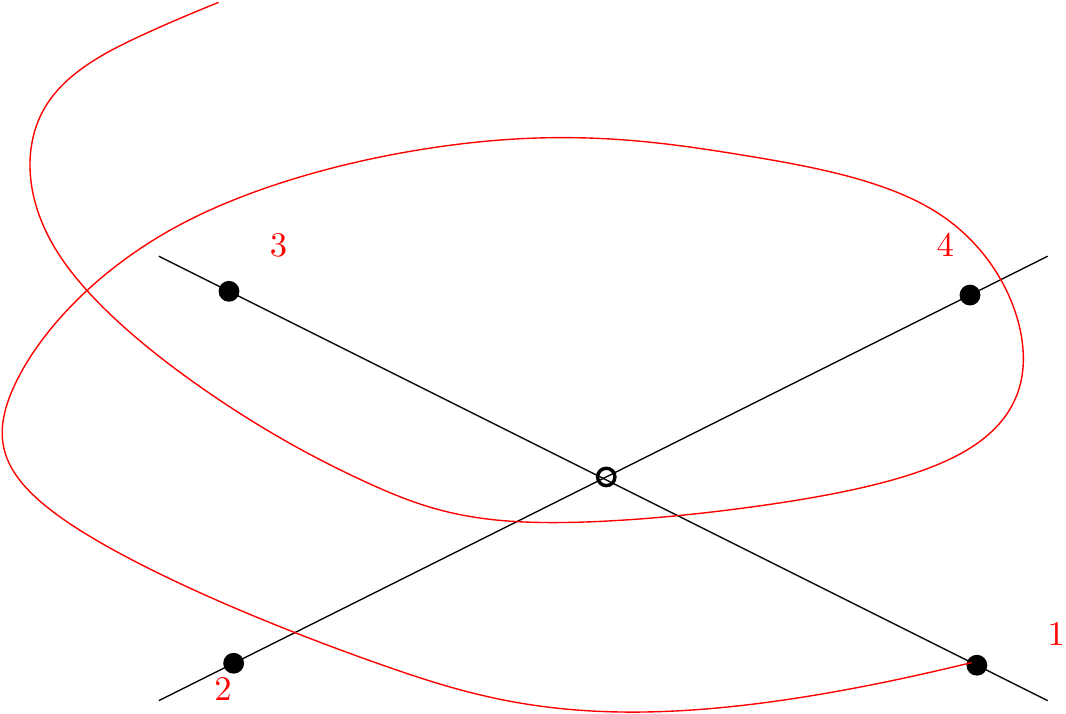}
\end{center}

\medskip

This gives a picture of an affine permutation $\tilde{s} = s \tau^{[a_1,a_2,a_3, a_4]}$, where $s(1)=3$ and $a_1=1$. This might be a little confusing, so we should take a moment to think about this. The integer $a_i$ is telling us to wrap completely around the circle $S_1$, $a_i$ times, then apply $s \in S_n$. If $a_i < 0$, then we move in the \emph{counterclockwise} direction (following our conventions), and if $a_i >0$, then we should wrap around \emph{clockwise} $a_i$ full rotations (through an angle of $a_i \cdot 2 \pi$). This of course can be thought of as, you guessed it, \textbf{complex multiplication}! In particular, we multiply the point we have labeled "$k = e^{2 \pi i k/n}$"
\[ e^{2 \pi i k/n} \mapsto e^{2 \pi i s(k)-k/n} e^{2 \pi i a_i} e^{2 \pi i k/n} \]
in terms of permutations of $\ZZ/n \ZZ$, we can think of this as
\[ k \otimes_{\CC} x^r \mapsto x^{a_i}  (k \otimes_{\CC} x^r) \]
which of course is an element of 
\[ \ZZ/n\ZZ \otimes_{\CC} \CC((x)) \subset S^1 \otimes_{\CC} \CC((x)).\]
Intuitively, as is suggested by the picture, one should think of this red loop around a cylinder corresponding to permutations $\tilde{s} \in \tilde{S}_n$, as "\emph{living over the vertex}" $x_j$ of a cellularly embedded graph in a Riemann surface $X$. Or, better, as a puncture or ramification point $x_j$ for some covering $p: X \to Y$. Moreover, we should think of $\sigma_x \in \tilde{S}_n$, the cyclic permutation of the monodromy group obtained from the constellations $C = [\sigma, \alpha, \phi]$, as being the lift of a cycle of $\sigma$ to an affine permutation via the map of 
\[ S_n \to \tilde{S}_n.\]

\begin{recall}
The cyclic permutation $\sigma_j$ was identified with the matrix
\[ \sigma_j= \begin{pmatrix}
	0 	  & 0 & 0 & \cdots & 0 & x_j \\
	1 	  & 0 & 0 & \cdots & 0 & 0 \\
	0 	  & 1 & 0 & \cdots & 0 & 0 \\
\vdots & \vdots & \vdots & \ddots & \vdots & \vdots \\
	0 	  & 0 & 0 &\cdots & 0 & 0 \\
	0 	  & 0 & 0 &\cdots & 1 & 0 
\end{pmatrix}_{n_i} \]
\end{recall}

Now we're cooking with gas. One more combinatorial/geometric/group theoretic realization of coverings of Riemann surfaces to add to our collection. Now, suppose we think of these loops given by elements $\tilde{s} \in \tilde{S}_n$ in the affine Weyl group. The Coxeter diagram of $\tilde{S}_n$ is

\begin{center}
\begin{tikzcd}
 & \circ_3 \arrow[r, no head] & \circ_4 \arrow[rd, no head] &  \\
\circ_2 \arrow[ru, no head] &  &  & \circ_5 \\
\circ_1 \arrow[u, no head] \arrow[rd, no head] &  &  &  \\
 & \circ_{n-1} \arrow[rruu, no head, dotted, bend right=49] &  & 
\end{tikzcd}
\end{center}

and as one often does, we can think of this as an unoriented version of the cyclic quiver with $n-1$ vertices corresponding to the vertex $x_j$ (where we have a puncture or ramification point) and the cyclic permutation $\sigma_j$. It is well known for the usual $A_n$ case for $\gl_n(\CC)$ with Weyl group $S_n \subset \GL_n(\CC)$, that the \textbf{indecomposable representations} of, say, the equioriented quiver\footnote{One can choose various orderings on the roots which corresponds to different orientations of the $A_n$ quiver, but this is not currently needed.}
\[ \bullet_1 \to \bullet_2 \to \cdots \to \bullet_{n-1} \]
correspond to the \textbf{positive roots} of the Lie algebra $\gl_n(\CC)$. A similar thing happens when we look at roots of $\gl_n(\FF)$, but we can also treat $\gl_n(\FF)$ as a $\CC$ algebra. In Section \ref{affine Schubert varieties}, the finite $\CC$-dimensional indecomposable representations will be identified with certain affine Schubert varieties following \cite{M}.

\section{Fibre-Products of Affine Grassmannians and Affine Flags}

In this section we will try to understand how the pullback of loop algebras gives a fibre product of affine Grassmannians and Affine Flag Varieties. So, Suppose we have a nontrivial morphism of compact Riemann surfaces (which is necessarily a covering)
\[ p: X \to Y \]
Let $\{x_1, x_2, ..., x_t\}$ be the ramification points. Now, to each $x_i$ we associate a disk 
\[ \DD_i = \DD_{x_i} = \Spec \CC((x_i)) = \Spec R_i \]
or a punctured disk if we are thinking of the $x_i$ as punctures,
\[ \DD_i^{\times} = \DD_{x_i}^{\times} = \Spec \CC((x_i)) = \Spec \FF_i \]
Now, let $d_i$ be the ramification index at $x_i$. Then we think of this as a vertex with $d_i$ half edges, and to it we associate the cyclic quiver $Q_i = Q_{x_i}$, with $d_i$ arrows (and vertices). Next we associate the matrix algebras
\[ M(\FF_i) := \Mat_{d_i \times d_i}(\FF_i) \supset \Mat_{d_i \times d_i}(R_i) = M(R_i) \]
and inside of these, we take the \textbf{Borel surface algebras} to be the Lie algebras of
\[ B_{-}(\FF_i) \supset B_{-}(R_i) \]
which are all of the form
\[ \fb((x_i))_{-} \supset \fb[[x_i]]_{-} \]
and are exactly the surface orders we constructed and are the completions of the hereditary path algebras in the normalization of the surface algebra. Next, using the construction of \emph{surface algebras and surface orders}, we take the pullback, and glue diagonal entries of the $M(\FF_i)$ using the data $C = [\sigma, \alpha, \phi]$. Now, over each $\DD_i$ we have the obvious bundles given by the matrix algebras $\Mat_{d_i \times d_i}(\FF_i)$. This gives us the following trivial bundle, with the various subbundles associated to the matrix subalgebras of $M(\FF_i)$ and $M(R_i)$.
\[ \mathcal{M}(X)  := \left( \prod_{i=1}^t M(\FF_i) \right) \otimes_{\CC} \cO(X) \]
where $\cO(X)$ is the structure sheaf of the curve $X$. Now, the pullback which we used to define the \emph{surface orders} defines a subbundle of the trivial bundle, 
The bundles given by the Borel subalgebras $\fb_i((x))$ give a subbundle of the trivial bundle given by their product
\[ \mathcal{B}(X)  := \left( \prod_{i=1}^t \fb((x_i)) \right) \otimes_{\CC} \cO(X) \]
We also get the corresponding quotient bundles. Now, when we glue the diagonals of the various matrix algebras over the various disks/punctured disks, we are essentially performing the construction explained in the section on the Gel'fand-Ponomarev Algebra, namely,
\[ \xymatrix{
\CC[[x,y]]/(x,y) \ar[r]  \ar[d] & \CC[[x]] \ar[d]^{g} \\
\CC[[y]] \ar[r]_{f} & \CC
}\]
In this construction, we are working over $\PP_{\CC}^1$, where we defined the "trivial dessin" to be the real segment $[0,1]$, and treated $\{0,1\}$ as punctures in order to get the fundamental group $\fF_2$, and the complete path algebra of the surface order 
\[ \CC \fF_2/\la xy, yx \ra = \CC \la \la x, y \ra \ra / \la xy, yx \ra \cong \CC[[x,y]]/(xy) \]
which we showed was isomorphic to the group algebra of the free group on two generators. 
Since we are working in $X = \PP_{\CC}^1-\{0,1\}$ with punctures, we took the completions, and so we are then choosing coordinates for the localization $\hO(X)_x$ and $\hO(X)_y$ in order to obtain isomorphisms $\hO(X)_x \cong \CC[[x]]$ and $\hO(X)_y \cong \CC[[y]]$. 

But, remember, this gluing is happening at \emph{every diagonal entry} of each matrix algebra over some puncture or ramification point. How should we interpret this? Remember, the surface algebras have nonzero cycles determined by the order of ramification at each $a_i$, and the size of the matrices for the surface orders were also determined in this way. and over a diagonal entry of one matrix (equivalently over a vertex of the quiver of the suface algebra) one glues to another diagonal entry (or vertex). What does this mean? Well, we should be thinking in terms of root systems at this point, i.e. in terms of the Cartan subalgebra 
\[ \fh[[x_i]] \subset \fb[[x_i]] \subset \fg[[x_i]] \subset \fg((x_i))\] 
in the hereditary order $\fb[[x_i]]$, in the normalization of the surface order, corresponding to that vertex. If we think of this in terms of quivers, this is gluing some idemptotents of the normalization of the surface algebra via some radical embedding sending each idempotent to a sum of two idempotents.

Now, gluing Cartan subalgebras $\fh[[x_i]] \subset \fh((x_i))$, means we are gluing root data for the corresponding Lie algebras $\fg[[x_i]]$ and the loop algebras $\fg((x_i))$. We can think of choosing an ordering of the roots as choosing the permutations $\sigma_{x_i}$ over each $x_i$. 

So in particular, each algebra $\fb[[x_i]] \subset \fb((x_i))$ identified to a vertex can be viewed as a product of affine flag varieties. Moreover, the representation theory of each individual loop algebra $\fg((x_i))$ and $\fg[[x_i]]$ is understood. Irreducible representations of $\GL_n(\CC)$ can be obtained using the Borel-Weil Theorem and taking sections of line bundles on the flag manifold. So at this point, we might want to try and do something similar for the affine loop algebras $\fg((x_i))$ and $\fg[[x_i]]$. This is good, but not enough. We are after more! We want the description of the representations of the pullback algebra given by the construction of the surface orders so that we can understand the representations of the affine Kac-Moody Lie algebras $\widehat{\fg}_i$
\[ \CC \mathbf{1} \to \widehat{\fg}_i \to \fg((x_i))  \]
Then, we would like to understand the pullback of all of the $\widehat{\fg}_i$.

\section{Lusztig's Isomorphism for the Gel'fand Ponomarev Algebra}

Lusztig gave the following very useful construction. Take $\Mat_{n \times n}(\CC)$ to be the matrix algebra with the action of $\GL_n(\CC)$ by conjugation. Lusztig gave an equivariant algebraic isomorphism of the variety $\mathcal{N} \subset \Mat_{n \times n}(\CC)$, of \textbf{nilpotent matrices}, and the \emph{opposite cell of a Schubert variety} in the \emph{affine Grassmannian} $\Gr(V) = \Gr(\FF^n)$. The basic idea is to take any $N \in \mathcal{N}$ to the "\emph{semi-infinite block-column matrix}"
\[ N \mapsto \begin{pmatrix}
N^{n-1} \\
\vdots \\
N^2 \\
N \\
I \\
0 \\
\vdots
\end{pmatrix}\]
where we have $\sum_{i=1}^{\infty} a_ie_i \in E_1$ represented as a semi-infinite column vector with $a_i \in R$, and $\Lambda = \Span_R \la v_1, ..., v_n \ra \subset \Gr(V)$ with $\{v_1, ..., v_n\}$ and $\FF$-basis of $V$. The $n \times n$ minors giving Pl\"{u}cker coordinates in $\Gr(V)$ restrict to polynomial functions on $\mathcal{N}$. This can be understood by defining the action of $N \in \cN$ on $V \cong \FF^n$ by
\[ \phi_N(v) = \frac{x^{n-1}}{1-x^{-1}N}(v) = \sum_{k=1}{n} t^{n-k}N^{k-1}(v) \]
and 
\[ \Phi: \mathcal{N} \to \Gr(V) \]
is given by 
\[ N \mapsto \begin{pmatrix}
N^{n-1} \\
\vdots \\
N^2 \\
N \\
I \\
0 \\
\vdots
\end{pmatrix}\]

For any $g \in \GL_n(\CC)$ we have
\[ \Phi g \cdot N) = \Phi(gNg^{-1}) = g\cdot \Phi(N)\]
is $\GL_n(\CC)$-equivarient.

\section{Lusztig's Isomorphism for Arbitrary Cyclic Quivers}

Lusztig later generalized this construction to a case which is of central importance to us. Let us take the cyclic quiver

\[
\tilde{A}_n := \vcenter{\hbox{  
		\begin{tikzpicture}[point/.style={shape=circle, fill=black, scale=.3pt,outer sep=3pt},>=latex]
		\node[point,label={above:$x_1$}] (0) at (0,0) {};
		\node[point,label={right:$x_2$}] (1) at (1.5,-.5) {};
		\node[point,label={left:$x_{|\sigma_i|}$}] (n) at (-1.5,-.5) {};
		\node[point,label={below:$x_{k}$}] (2) at (0,-4) {};
		\node[point,label={right:$x_{k-1}$}] (3) at (1.5,-3.5) {};
		\node[point,label={left:$x_{k+1}$}] (4) at (-1.5,-3.5) {};
		
		% \draw[dotted] (2.2,-1.65)--(2.2,-2.4);
		% \draw[dotted] (-2.2,-1.65)--(-2.2,-2.4);
		
		\path[->]
		(4) edge [dashed,bend left=45]  (n)
		(1) edge [dashed,bend left=45]  (3)
		(0) edge [bend left=15] node[midway, above] {$a_1$} (1)
		(3) edge [bend left=15] node[midway, below] {$a_{k}$} (2)
		(2) edge [bend left=15] node[midway, below] {$a_{k+1}$} (4)
		(n) edge [bend left=15] node[midway, above] {$a_{|\sigma_i|}$} (0);
		\end{tikzpicture} 
}}
\]
with corresponding path algebra $k\tilde{A}_{n}$. Define the \textbf{dimension vector} $\bd = (d_1, ..., d_n)$ to be an assignment of nonnegative integers to the vertices. A $\bd$-dimensional representation of this quiver is then an assignment of linear maps from
\[ (A_1, ..., A_n) \in \prod_{i=1}^{n} \Mat_{d_{i-1} \times d_i}(\CC) \]
to each arrow with indices taken modulo $n$. This is of course an affine space. Then there is an action of $\GL(\bd) = \prod_{i=1}^n \GL_{d_i}(\CC)$ by conjugation of the matrices, i.e. by base change on the $d_i$ dimensional $\CC$-vector space associated to each vertex. The space of \textbf{nilpotent representation} is
\[ \mathcal{A}(\bd) = \{(A_1, ..., A_n) : \ \prod_{i=1}^n M_{n-i} = 0\}.\]
We can take any cyclic permutation of this as well. Now define a map
\[ (A_1, ..., A_n) \mapsto \begin{pmatrix}
0 & 0 & 0 & \cdots & 0 & 0  \\
A_1 & 0 & 0 & \cdots & 0 & 0 \\
0 & A_2 & 0 & \cdots & 0 &  \\
\vdots & \vdots & \vdots & \ddots & \vdots & \vdots \\
0 & 0 & 0 & \cdots & 0 & 0 \\
x^{-1} \cdot A_n & 0 & 0 & \cdots & A_{n-1} & 0
\end{pmatrix} \in \GL_n(\FF). \]
Next, define $R = \CC[[x]]$-lattices
\[ \Lambda_1 = \begin{pmatrix}
A_1^{[kn-n]} & 0 & \cdots & 0 \\
\vdots & \vdots & & \vdots \\
A_nA_1 & 0 & \cdots & 0 \\
A_1 & 0 & \cdots & 0 \\
I_1 & 0 & \cdots & 0 \\
0 & I_2 & \cdots & 0 \\
\vdots & \vdots & \ddots & \vdots \\
0 & 0 & \cdots & I_n \\
0 & 0 & \cdots & 0 \\
\vdots & \vdots & & \vdots 
\end{pmatrix}, \quad \Lambda_2 = \begin{pmatrix}
A_2^{[kn-n+1]} & 0 & \cdots & 0 \\
\vdots & \vdots & & \vdots \\
A_nA_1A_2 & 0 & \cdots & 0 \\
A_1A_2 & 0 & \cdots & 0 \\
A_2 & 0 & \cdots & 0 \\
I_2 & 0 & \cdots & 0 \\
0 & I_3 & \cdots & 0 \\
\vdots & \vdots & \ddots & \vdots \\
0 & 0 & \cdots & I_n \\
0 & 0 & \cdots & 0 \\
\vdots & \vdots & & \vdots 
\end{pmatrix}, \quad \cdots \]
where in Magyar's notation we have $A_{j+nk} = A_j$ and $A_j^{[k]}:= M_{j-k+1} \cdots M_{j-1}M_j$. In \cite{M} it is then shown that the image of $\Phi$ embeds $\mathcal{A}(\bd)$ into a partial flag variety which he calls $\Fl(\bd, V)$. In our notation, we may realize this as the partial flag in $\fg((x_i))$ as follows. Choose any composition of $n_i$, the order of the ramification at $x_i$ in the surface $X$ (i.e. the number of arrows in the cyclic quiver over $x_i$ and the dimension of the matrix algebra $\fg((x_i))$). Let us call this composition
\[ \bd = d_1+d_2+ \cdots d_k \]
Next, define a flag in $\fg((x_i))$ in the usual way by block upper triangular matrices with blocks given by the composition $\bd$. 

Now, once we have a partial flag $\Fl(V_i)$ for each $x_i$, where $\dim_{\FF}(V_i) = n_i = \sum_{j=1}^{k_i} d_j$, then we may construct a trivial bundle
\[ \cF \otimes_{\CC} \cO(X) \]
where $\cF = \prod_{x_i} \Fl(V_i)$. Embedded in this trivial bundle is the trivial trivial affine Schubert bundle, which is given by the product of the images of the maps $\Phi_i: \mathcal{A}_i(\bd_i) \to \Fl(V_i)$ giving the Schubert varieties for the cyclic quivers over the $x_i$. Now, the gluing of the diagonal entries of each $\fg((x_i))$ via the surface order pullback construction induces a pullback on these Schubert varieties as well. Moreover, for each partial flag $\Fl(V_i) \cong \GL_{n_i}(\FF_i)/B(\FF_i)_{+}$, we may define a pullback of homogeneous spaces under this construction as well.

\section{Next Steps}

Now that we have a meaningful way of understanding representations of pullbacks of loop groups and loop algebras in terms of geometry via the Schubert varieties construction, and we have a very clear connection to the Galois groups of extensions $\cG(\CC(X)/\CC(Y))$ for maps of function fields $\CC(Y) \to \CC(X)$ coming from coverings of Riemann surfaces $p:X \to Y$, we would like to begin understanding this geometry a little better and understanding what this representation theory-geometry correspondence tells us. It is too much to set up the theory of representation varieties of quivers and their corresponding moduli spaces here, so in the next paper, we will describe the indecomposable representations of surface algebras, the representation varieties, and the corresponding representation theory of surface orders via Magyar's setup of Lusztig's isomorphism. We will also give a nice description of their moduli spaces.

\clearpage
\textcolor{darkred}{\hrule}

{00}

\end{document}